\documentclass[thmsa,onecolumn,10pt]{article}
\usepackage{amssymb}
\usepackage{amsfonts}
\usepackage{graphicx}
\usepackage{amsmath}

\newtheorem{theorem}{Theorem}

\newtheorem{corollary}[theorem]{Corollary}

\newtheorem{lemma}[theorem]{Lemma}

\newtheorem{proposition}[theorem]{Proposition}

\newenvironment{proof}[1][Proof]{\textbf{#1.} }{\ \rule{0.5em}{0.5em}}

\begin{document}

\title{Identities in law between quadratic functionals of bivariate Gaussian
processes, through Fubini theorems and symmetric projections}
\author{Giovanni PECCATI\thanks{%
Laboratoire de Statistique Th\'{e}orique et Appliqu\'{e}e, Universit\'{e}
Paris VI, France. E-mail: \texttt{giovanni.peccati@gmail.com}}, and Marc YOR%
\thanks{%
Laboratoire de Probabilit\'{e}s et Mod\`{e}les Al\'{e}atoires, Universit\'{e}%
s Paris VI and Paris VII, France and Institut Universitaire de France.}}
\date{January 27, 2005}
\maketitle

\begin{abstract}
We present three new identities in law for quadratic functionals of
conditioned bivariate Gaussian processes. In particular, our results provide
a two-parameter generalization of a celebrated identity in law, involving
the path variance of a Brownian bridge, due to Watson (1961). The proof is
based on ideas from a recent note by J. R. Pycke (2005) and on the
stochastic Fubini theorem for general Gaussian measures proved in Deheuvels
\textit{et al. }(2004).


\textbf{AMS\ 2000 classification : }60515, 60E10.
\end{abstract}

\section{Introduction}

Let $b\left( s\right) $, $s\in \left[ 0,1\right] $, be a standard Brownian
bridge on $\left[ 0,1\right] $, from $0$ to $0$, and let $b_{1}$ and $b_{2}$
be two independent copies of $b$. The aim of this note is to prove several
bivariate generalizations of the following identity in law for the path
variance of $b$,
\begin{equation}
\int_{0}^{1}\left( b\left( s\right) -\int_{0}^{1}b\left( u\right) du\right)
^{2}ds\overset{law}{=}\frac{1}{4}\int_{0}^{1}\left[ b_{1}\left( s\right)
^{2}+b_{2}\left( s\right) ^{2}\right] ds,  \label{WI}
\end{equation}%
known as \textit{Watson's (duplication) identity }(see \cite{Watson}; the
reader is also referred to \cite{ShiYor} for a detailed probabilistic
discussion of (\ref{WI})). More specifically, our aim is to establish a
result analogous to (\ref{WI}) for the path variance of a \textit{bivariate
tied-down Brownian bridge} $\mathbf{B}_{0}$ on $\left[ 0,1\right] ^{2}$,
i.e. a process having the law of a standard Brownian sheet $\mathbf{W}$
conditioned to vanish on the edges of the square $\left[ 0,1\right] ^{2}$.
As discussed below, our bivariate generalizations of (\ref{WI}) involve four
different types of \textquotedblleft bridges\textquotedblright\ naturally
attached to a given Brownian sheet $\mathbf{W}$. These four processes, along
with the laws of their quadratic functionals, have been recently studied in
\cite{DPY}.

The relatively simple proof of our main result uses extensively the general
\textit{stochastic Fubini theorem}, for quadratic functionals of Gaussian
processes, proved in \cite{DPY} (but see also \cite{DMartinYor}), and has
been inspired by the recent proof of Watson's identity given in \cite{Pycke}%
. Such a proof is mainly based on a decomposition of the path of\ the random
function $t\mapsto b\left( t\right) $ into the orthogonal sum of its
symmetric and antisymmetric parts, around the pole $t=1/2$. We will see how
this kind of decomposition can be naturally extended to the framework of
bivariate functions.

The present note is organized as follows. In Section 2 we introduce some
notation. In Section 3, we state a version of the stochastic Fubini Theorem
which is well adapted to the framework of this paper and we provide an
alternative proof of such a result, based on the calculation of cumulants
for double Wiener integrals. In Section 4 the main Theorem is stated and
proved. Eventually, in Section 5 we apply our results to calculate: (a) the
explicit Laplace transform of some quadratic functionals of bivariate
Gaussian processes, and (b) the explicit Fourier transform of some double
stochastic integrals with respect to conditioned bivariate processes. This
completes part of the results obtained in \cite{DPY} and \cite{JuliaNualart}.

\section{General notation}

For the rest of the paper, we will study Gaussian processes that can be
expressed as suitable transformations of a standard Brownian motion or of a
standard Brownian sheet. In particular, we will adopt the following notation:

\begin{description}
\item[ ] -- $W=\left\{ W\left( t\right) :t\in \left[ 0,1\right] \right\} $
is a standard Brownian motion on $\left[ 0,1\right] $, initialized at 0;

\item[ ] -- $b=\left\{ b\left( t\right) :t\in \left[ 0,1\right] \right\} $
is a standard Brownian bridge on $\left[ 0,1\right] $, from $0$ to $0$;

\item[ ] -- $\mathbf{W=}\left\{ \mathbf{W}\left( t_{1},t_{2}\right) :\left(
t_{1},t_{2}\right) \in \left[ 0,1\right] ^{2}\right\} $ is a standard
Brownian sheet on $\left[ 0,1\right] ^{2}$ vanishing on the axes, that is, $%
\mathbf{W}$ is a centered Gaussian process such that, for every $\left(
t_{1},t_{2}\right) $, $\left( s_{1},s_{2}\right) \in \left[ 0,1\right] ^{2}$%
\begin{equation*}
\mathbb{E}\left[ \mathbf{W}\left( t_{1},t_{2}\right) \mathbf{W}\left(
s_{1},s_{2}\right) \right] =\left( t_{1}\wedge s_{1}\right) \times \left(
t_{2}\wedge s_{2}\right) ;
\end{equation*}

\item[ ] -- $\mathbf{B}^{\left( \mathbf{W}\right) }=\left\{ \mathbf{B}%
^{\left( \mathbf{W}\right) }\left( t_{1},t_{2}\right) :\left(
t_{1},t_{2}\right) \in \left[ 0,1\right] ^{2}\right\} $ is the canonical
\textit{bivariate Brownian bridge} associated to $\mathbf{W}$, i.e.%
\begin{equation*}
\mathbf{B}^{\left( \mathbf{W}\right) }\left( t_{1},t_{2}\right) =\mathbf{W}%
\left( t_{1},t_{2}\right) -t_{1}t_{2}\mathbf{W}\left( 1,1\right) ;
\end{equation*}

\item[ ] -- $\mathbf{B}_{0}^{\left( \mathbf{W}\right) }=\left\{ \mathbf{B}%
_{0}^{\left( \mathbf{W}\right) }\left( t_{1},t_{2}\right) :\left(
t_{1},t_{2}\right) \in \left[ 0,1\right] ^{2}\right\} $ is the canonical
\textit{bivariate tied down Brownian bridge} associated to $\mathbf{W}$, i.e.%
\begin{equation*}
\mathbf{B}_{0}^{\left( \mathbf{W}\right) }\left( t_{1},t_{2}\right) =\mathbf{%
W}\left( t_{1},t_{2}\right) -t_{1}\mathbf{W}\left( 1,t_{2}\right) -t_{2}%
\mathbf{W}\left( t_{1},1\right) +t_{1}t_{2}\mathbf{W}\left( 1,1\right) ;
\end{equation*}

\item[ ] -- $\mathbf{K}^{\left( \mathbf{W},i\right) }=\left\{ \mathbf{K}%
^{\left( \mathbf{W},i\right) }\left( t_{1},t_{2}\right) :\left(
t_{1},t_{2}\right) \in \left[ 0,1\right] ^{2}\right\} $, $i=1,2$, are the
two canonical \textit{Kiefer fields} (or \textit{asymmetric bivariate bridges%
}) associated to $\mathbf{W}$, i.e.%
\begin{eqnarray*}
\mathbf{K}^{\left( \mathbf{W},1\right) }\left( t_{1},t_{2}\right) &=&\mathbf{%
W}\left( t_{1},t_{2}\right) -t_{1}\mathbf{W}\left( 1,t_{2}\right) \\
\mathbf{K}^{\left( \mathbf{W},2\right) }\left( t_{1},t_{2}\right) &=&\mathbf{%
W}\left( t_{1},t_{2}\right) -t_{2}\mathbf{W}\left( t_{1},1\right) .
\end{eqnarray*}
\end{description}

\bigskip

\bigskip We assume that all the previous objects are defined on the same
probability space $\left( \Omega ,\mathcal{F},\mathbb{P}\right) $.

\begin{description}
\item[ ] \textbf{Remarks -- }(i) Conditionally on the event $\left\{ \mathbf{%
W}\left( 1,1\right) =0\right\} $, $\mathbf{W}$ is distributed as the
unconditioned process $\mathbf{B}^{\left( \mathbf{W}\right) }.$ Moreover,
for every $\left( t_{1},t_{2}\right) $, $\left( s_{1},s_{2}\right) \in \left[
0,1\right] ^{2}$,%
\begin{equation}
\mathbb{E}\left[ \mathbf{B}^{\left( \mathbf{W}\right) }\left(
t_{1},t_{2}\right) \mathbf{B}^{\left( \mathbf{W}\right) }\left(
s_{1},s_{2}\right) \right] =\left( t_{1}\wedge s_{1}\right) \times \left(
t_{2}\wedge s_{2}\right) -t_{1}s_{1}t_{2}s_{2}.  \label{prod0}
\end{equation}

\item[ ] (ii) Conditionally on the event $\left\{ \mathbf{W}\left(
1,t\right) =\mathbf{W}\left( t,1\right) =0,\text{ \ }\forall t\in \left[ 0,1%
\right] \right\} $, $\mathbf{W}$ is distributed as the unconditioned process
$\mathbf{B}_{0}^{\left( \mathbf{W}\right) }$. In particular, for $\left(
t_{1},t_{2}\right) $, $\left( s_{1},s_{2}\right) \in \left[ 0,1\right] ^{2}$%
,
\begin{eqnarray}
\mathbb{E}\left[ \mathbf{B}_{0}^{\left( \mathbf{W}\right) }\left(
t_{1},t_{2}\right) \mathbf{B}_{0}^{\left( \mathbf{W}\right) }\left(
s_{1},s_{2}\right) \right]  &=&\mathbb{E}\left[ b\left( t_{1}\right) b\left(
s_{1}\right) \right] \times \mathbb{E}\left[ b\left( t_{2}\right) b\left(
s_{2}\right) \right]   \label{prod1} \\
&=&\left( t_{1}\wedge s_{1}-t_{1}s_{1}\right) \times \left( t_{2}\wedge
s_{2}-t_{2}s_{2}\right) .  \notag
\end{eqnarray}

\item[ ] (iii) Conditionally on $\left\{ \mathbf{W}\left( 1,t\right) =0,%
\text{ \ }\forall t\in \left[ 0,1\right] \right\} $, $\mathbf{W}$ is
distributed as the unconditioned process $\mathbf{K}^{\left( \mathbf{W}%
,1\right) }$, and moreover, for $\left( t_{1},t_{2}\right) $, $\left(
s_{1},s_{2}\right) \in \left[ 0,1\right] ^{2}$,
\begin{eqnarray}
\mathbb{E}\left[ \mathbf{K}^{\left( \mathbf{W},1\right) }\left(
t_{1},t_{2}\right) \mathbf{K}^{\left( \mathbf{W},1\right) }\left(
s_{1},s_{2}\right) \right]  &=&\mathbb{E}\left[ b\left( t_{1}\right) b\left(
s_{1}\right) \right] \times \mathbb{E}\left[ W\left( t_{2}\right) W\left(
s_{2}\right) \right]   \label{prod2} \\
&=&\left( t_{1}\wedge s_{1}-t_{1}s_{1}\right) \times \left( t_{2}\wedge
s_{2}\right) .  \notag
\end{eqnarray}

\item[ ] (iv) Conditionally on $\left\{ \mathbf{W}\left( t,1\right) =0,\text{
\ }\forall t\in \left[ 0,1\right] \right\} $, $\mathbf{W}$ is distributed as
the unconditioned process $\mathbf{K}^{\left( \mathbf{W},2\right) }$, and
moreover, for $\left( t_{1},t_{2}\right) $, $\left( s_{1},s_{2}\right) \in %
\left[ 0,1\right] ^{2}$,
\begin{eqnarray}
\mathbb{E}\left[ \mathbf{K}^{\left( \mathbf{W},2\right) }\left(
t_{1},t_{2}\right) \mathbf{K}^{\left( \mathbf{W},2\right) }\left(
s_{1},s_{2}\right) \right]  &=&\mathbb{E}\left[ W\left( t_{1}\right) W\left(
s_{1}\right) \right] \times \mathbb{E}\left[ b\left( t_{2}\right) b\left(
s_{2}\right) \right]   \label{prod3} \\
&=&\left( t_{1}\wedge s_{1}\right) \times \left( t_{2}\wedge
s_{2}-t_{2}s_{2}\right) .  \notag
\end{eqnarray}
\end{description}

\section{Stochastic Fubini Theorems}

The following stochastic Fubini theorem will be useful for the proof of our
main results. As shown in \cite{DMartinYor} and \cite{DPY}, stochastic
Fubini theorems for general Gaussian measures can be easily proved by means
of a Laplace transform argument. Here, we shall present an alternative
proof, which is based on the so called \textit{diagram formulae} (see e.g.
\cite{Sur}) for the cumulants of double Wiener integrals. Note that, in what
follows, we will write $d\lambda ^{m}$, $m\geq 1$, to indicate Lebesgue
measure on $\Re ^{m}$.

\begin{theorem}[Stochastic Fubini Theorem]
Under the above assumptions and notation, for every $\phi \in L^{2}\left( %
\left[ 0,1\right] ^{4},d\lambda ^{4}\right) $ there exist two measurable
random functions
\begin{equation}
\left\{ \int_{\left[ 0,1\right] ^{2}}\phi \left(
t_{1},t_{2},x_{1},x_{2}\right) \mathbf{W}\left( dx_{1},dx_{2}\right) :\left(
t_{1},t_{2}\right) \in \left[ 0,1\right] ^{2}\right\}  \label{R1}
\end{equation}%
and
\begin{equation}
\left\{ \int_{\left[ 0,1\right] ^{2}}\phi \left(
x_{1},x_{2},t_{1},t_{2}\right) \mathbf{W}\left( dx_{1},dx_{2}\right) :\left(
t_{1},t_{2}\right) \in \left[ 0,1\right] ^{2}\right\} .  \label{R2}
\end{equation}%
Moreover, the following distributional identity holds
\begin{eqnarray}
&&\int_{\left[ 0,1\right] ^{2}}\left[ \int_{\left[ 0,1\right] ^{2}}\phi
\left( t_{1},t_{2},x_{1},x_{2}\right) \mathbf{W}\left( dx_{1},dx_{2}\right) %
\right] ^{2}dt_{1}dt_{2}  \label{SFT} \\
&&\overset{law}{=}\int_{\left[ 0,1\right] ^{2}}\left[ \int_{\left[ 0,1\right]
^{2}}\phi \left( x_{1},x_{2},t_{1},t_{2}\right) \mathbf{W}\left(
dx_{1},dx_{2}\right) \right] ^{2}dt_{1}dt_{2}.  \notag
\end{eqnarray}
\end{theorem}

\begin{proof}
The existence of the two measurable random functions (\ref{R1}) and (\ref{R2}%
) follows from standard arguments. To obtain (\ref{SFT}), start by defining
the two kernels (contractions) on $\left[ 0,1\right] ^{4}$%
\begin{eqnarray*}
\Phi _{1}\left( t_{1},t_{2};s_{1},s_{2}\right)  &=&\int_{\left[ 0,1\right]
^{2}}dx_{1}dx_{2}\phi \left( x_{1},x_{2},t_{1},t_{2}\right) \phi \left(
x_{1},x_{2},s_{1},s_{2}\right)  \\
\Phi _{2}\left( t_{1},t_{2};s_{1},s_{2}\right)  &=&\int_{\left[ 0,1\right]
^{2}}dx_{1}dx_{2}\phi \left( t_{1},t_{2},x_{1},x_{2}\right) \phi \left(
s_{1},s_{2},x_{1},x_{2}\right) .
\end{eqnarray*}%
Then, a simple application of the multiplication formula for Wiener
integrals (see for instance \cite[p. 211]{DMM}) shows that%
\begin{eqnarray}
\int_{\left[ 0,1\right] ^{2}}\left[ \int_{\left[ 0,1\right] ^{2}}\phi \left(
t_{1},t_{2},x_{1},x_{2}\right) \mathbf{W}\left( dx_{1},dx_{2}\right) \right]
^{2}dt_{1}dt_{2} &=&\left\Vert \phi \right\Vert ^{2}+I_{2}^{\mathbf{W}%
}\left( \Phi _{1}\right)   \label{Wre} \\
\int_{\left[ 0,1\right] ^{2}}\left[ \int_{\left[ 0,1\right] ^{2}}\phi \left(
x_{1},x_{2},t_{1},t_{2}\right) \mathbf{W}\left( dx_{1},dx_{2}\right) \right]
^{2}dt_{1}dt_{2} &=&\left\Vert \phi \right\Vert ^{2}+I_{2}^{\mathbf{W}%
}\left( \Phi _{2}\right) ,  \notag
\end{eqnarray}%
where $I_{2}^{\mathbf{W}}\left( \cdot \right) $ stands for a standard double
Wiener integral with respect to $\mathbf{W}$ (see again \cite{DMM}). Now
define $\chi _{m}\left( Y\right) $, $m\geq 1$, to be the $m$-th cumulant of
a given real valued random variable $Y$ (see e.g. \cite{Sur}). We recall
that the law of a double Wiener integral is determined by its cumulants (see
e.g. \cite{Slud}). Moreover, we can apply the well known \textit{diagram
formulae} for cumulants of multiple stochastic integrals (as presented, for
instance, in \cite[Proposition 9 and Corollary 1]{Rota Wall}) to obtain that
for every $m\geq 2$ there exists a combinatorial coefficient $c_{m}>0$
(independent of $\phi $) such that
\begin{eqnarray}
\chi _{m}\left( I_{2}^{\mathbf{W}}\left( \Phi _{1}\right) \right)
&=&c_{m}\int_{\left[ 0,1\right] ^{2m}}\left( d\lambda ^{2}\right) ^{\otimes
m}\text{ }\Phi _{1}\left( x_{1}^{\left( 1\right) },x_{2}^{\left( 1\right)
};x_{1}^{\left( 2\right) },x_{2}^{\left( 2\right) }\right)   \label{cumu} \\
&&\times \Phi _{1}\left( x_{1}^{\left( 2\right) },x_{2}^{\left( 2\right)
};x_{1}^{\left( 3\right) },x_{2}^{\left( 3\right) }\right) \times \cdot
\cdot \cdot \times \Phi _{1}\left( x_{1}^{\left( m\right) },x_{2}^{\left(
m\right) };x_{1}^{\left( 1\right) },x_{2}^{\left( 1\right) }\right)   \notag
\\
&=&c_{m}\int_{\left[ 0,1\right] ^{2m}}\left( d\lambda ^{2}\right) ^{\otimes
m}\text{ }\Phi _{2}\left( x_{1}^{\left( 1\right) },x_{2}^{\left( 1\right)
};x_{1}^{\left( 2\right) },x_{2}^{\left( 2\right) }\right)   \notag \\
&&\times \Phi _{2}\left( x_{1}^{\left( 2\right) },x_{2}^{\left( 2\right)
};x_{1}^{\left( 3\right) },x_{2}^{\left( 3\right) }\right) \times \cdot
\cdot \cdot \times \Phi _{2}\left( x_{1}^{\left( m\right) },x_{2}^{\left(
m\right) };x_{1}^{\left( 1\right) },x_{2}^{\left( 1\right) }\right)   \notag
\\
&=&\chi _{m}\left( I_{2}^{\mathbf{W}}\left( \Phi _{2}\right) \right) ,
\notag
\end{eqnarray}%
where the second equality can be proved by using a standard (deterministic)
Fubini theorem, as well as the definition of $\Phi _{1}$ and $\Phi _{2}$.
Since (\ref{cumu}) holds for every $m$, we obtain that $I_{2}^{\mathbf{W}%
}\left( \Phi _{1}\right) \overset{law}{=}I_{2}^{\mathbf{W}}\left( \Phi
_{2}\right) $, and the proof of (\ref{SFT}) is therefore concluded, due to (%
\ref{Wre}).
\end{proof}

\bigskip

As shown in \cite{DPY}, by specializing (\ref{SFT}) to the kernels%
\begin{eqnarray*}
\phi ^{\left( 1\right) }\left( t_{1},t_{2};x_{1},x_{2}\right) &=&\mathbf{1}_{%
\left[ 0,t_{1}\right] }\left( x_{1}\right) \mathbf{1}_{\left[ 0,t_{2}\right]
}\left( x_{2}\right) -t_{1}t_{2} \\
\phi ^{\left( 2\right) }\left( t_{1},t_{2};x_{1},x_{2}\right) &=&\mathbf{1}_{%
\left[ 0,t_{1}\right] }\left( x_{1}\right) \mathbf{1}_{\left[ 0,t_{2}\right]
}\left( x_{2}\right) -t_{1}\mathbf{1}_{\left[ 0,t_{2}\right] }\left(
x_{2}\right) -t_{2}\mathbf{1}_{\left[ 0,t_{1}\right] }\left( x_{1}\right)
+t_{1}t_{2} \\
\phi ^{\left( 3\right) }\left( t_{1},t_{2};x_{1},x_{2}\right) &=&\mathbf{1}_{%
\left[ 0,t_{1}\right] }\left( x_{1}\right) \mathbf{1}_{\left[ 0,t_{2}\right]
}\left( x_{2}\right) -t_{1}\mathbf{1}_{\left[ 0,t_{2}\right] }\left(
x_{2}\right) \\
\phi ^{\left( 4\right) }\left( t_{1},t_{2};x_{1},x_{2}\right) &=&\mathbf{1}_{%
\left[ 0,t_{1}\right] }\left( x_{1}\right) \mathbf{1}_{\left[ 0,t_{2}\right]
}\left( x_{2}\right) -t_{2}\mathbf{1}_{\left[ 0,t_{1}\right] }\left(
x_{1}\right) ,
\end{eqnarray*}%
we obtain the following

\begin{corollary}
Let the above notation and assumptions prevail. Then,%
\begin{equation}
\int_{\left[ 0,1\right] ^{2}}\mathbf{B}^{\left( \mathbf{W}\right) }\left(
t_{1},t_{2}\right) ^{2}dt_{1}dt_{2}\overset{law}{=}\int_{\left[ 0,1\right]
^{2}}\left[ \mathbf{W}\left( t_{1},t_{2}\right) -\int_{\left[ 0,1\right]
^{2}}\mathbf{W}\left( u_{1},u_{2}\right) du_{1}du_{2}\right] ^{2}dt_{1}dt_{2}
\label{Fub1}
\end{equation}

\begin{gather}
\int_{\left[ 0,1\right] ^{2}}\mathbf{B}_{0}^{\left( \mathbf{W}\right)
}\left( t_{1},t_{2}\right) ^{2}dt_{1}dt_{2}\overset{law}{=}\int_{\left[ 0,1%
\right] ^{2}}\left[ \mathbf{W}\left( t_{1},t_{2}\right) -\int_{\left[ 0,1%
\right] }\mathbf{W}\left( t_{1},u_{2}\right) du_{2}\right.  \label{Fub2} \\
\text{ \ \ \ \ \ \ \ \ \ \ \ \ \ \ \ \ \ \ \ \ \ \ \ \ \ \ \ \ \ \ \ \ }%
\left. -\int_{\left[ 0,1\right] }\mathbf{W}\left( u_{1},t_{2}\right)
du_{1}+\int_{\left[ 0,1\right] ^{2}}\mathbf{W}\left( u_{1},u_{2}\right)
du_{1}du_{2}\right] ^{2}dt_{1}dt_{2}  \notag
\end{gather}

\begin{equation}
\int_{\left[ 0,1\right] ^{2}}\mathbf{K}^{\left( \mathbf{W},1\right) }\left(
t_{1},t_{2}\right) ^{2}dt_{1}dt_{2}\overset{law}{=}\int_{\left[ 0,1\right]
^{2}}\left[ \mathbf{W}\left( t_{1},t_{2}\right) -\int_{\left[ 0,1\right] }%
\mathbf{W}\left( u_{1},t_{2}\right) du_{1}\right] ^{2}dt_{1}dt_{2}
\label{Fubi3}
\end{equation}%
\begin{equation}
\int_{\left[ 0,1\right] ^{2}}\mathbf{K}^{\left( \mathbf{W},2\right) }\left(
t_{1},t_{2}\right) ^{2}dt_{1}dt_{2}\overset{law}{=}\int_{\left[ 0,1\right]
^{2}}\left[ \mathbf{W}\left( t_{1},t_{2}\right) -\int_{\left[ 0,1\right] }%
\mathbf{W}\left( t_{1},u_{2}\right) du_{2}\right] ^{2}dt_{1}dt_{2}
\label{Fub4}
\end{equation}
\end{corollary}

\section{Bivariate Watson's identities}

\subsection{Main results}

The next Theorem, which contains the announced bivariate versions of
Watson's duplication identity (\ref{WI}), is the main result of the section.
Note that each of the three parts of the statement involves a different
notion of path variance for the process $\mathbf{B}_{0}^{\left( \mathbf{W}%
\right) }$.

\begin{theorem}
Let $\mathbf{W}$ be a standard Brownian sheet on $\left[ 0,1\right] ^{2}$,
and let $\mathbf{W}_{i}$, $i=1,2,3,4$, be four independent copies of $%
\mathbf{W}.$ Then,

\begin{enumerate}
\item
\begin{align*}
& \int_{\left[ 0,1\right] ^{2}}\left[ \mathbf{B}_{0}^{\left( \mathbf{W}%
\right) }\left( t_{1},t_{2}\right) -\int_{\left[ 0,1\right] ^{2}}\mathbf{B}%
_{0}^{\left( \mathbf{W}\right) }\left( u_{1},u_{2}\right) du_{1}du_{2}\right]
^{2}dt_{1}dt_{2} \\
& \overset{\text{law}}{=}\frac{1}{16}\int_{\left[ 0,1\right] ^{2}}\left[
\mathbf{B}^{\left( \mathbf{W}_{1}\right) }\left( t_{1},t_{2}\right) ^{2}+%
\mathbf{K}^{\left( \mathbf{W}_{2},1\right) }\left( t_{1},t_{2}\right) ^{2}+%
\mathbf{K}^{\left( \mathbf{W}_{3},2\right) }\left( t_{1},t_{2}\right) ^{2}+%
\mathbf{B}_{0}^{\left( \mathbf{W}_{4}\right) }\left( t_{1},t_{2}\right) ^{2}%
\right] dt_{1}dt_{2}
\end{align*}

\item
\begin{eqnarray*}
&&\int_{\left[ 0,1\right] ^{2}}\left[ \mathbf{B}_{0}^{\left( \mathbf{W}%
\right) }\left( t_{1},t_{2}\right) -\int_{0}^{1}\mathbf{B}_{0}^{\left(
\mathbf{W}\right) }\left( t_{1},u_{2}\right) du_{2}\right] ^{2}dt_{1}dt_{2}
\\
&&\text{ \ \ \ \ \ \ \ \ \ \ \ \ \ \ \ \ \ \ \ \ \ \ \ \ \ \ \ \ \ \ \ \ \ \
\ \ \ \ }\overset{\text{law}}{=}\frac{1}{4}\int_{\left[ 0,1\right] ^{2}}%
\left[ \mathbf{B}_{0}^{\left( \mathbf{W}_{1}\right) }\left(
t_{1},t_{2}\right) ^{2}+\mathbf{B}_{0}^{\left( \mathbf{W}_{2}\right) }\left(
t_{1},t_{2}\right) ^{2}\right] dt_{1}dt_{2}
\end{eqnarray*}

\item
\begin{eqnarray*}
&&\int_{\left[ 0,1\right] ^{2}}\left[ \mathbf{B}_{0}^{\left( \mathbf{W}%
\right) }\left( t_{1},t_{2}\right) -\int_{0}^{1}\mathbf{B}_{0}^{\left(
\mathbf{W}\right) }\left( t_{1},u_{2}\right) du_{2}\right. \\
&&\text{ \ }\left. -\int_{0}^{1}\mathbf{B}_{0}^{\left( \mathbf{W}\right)
}\left( u_{1},t_{2}\right) du_{1}+\int_{\left[ 0,1\right] ^{2}}\mathbf{B}%
_{0}^{\left( \mathbf{W}\right) }\left( u_{1},u_{2}\right) du_{1}du_{2}\right]
^{2}dt_{1}dt_{2} \\
&&\text{ \ \ \ \ \ \ \ \ \ \ \ \ }\overset{law}{=}\frac{1}{16}\int_{\left[
0,1\right] ^{2}}\sum_{i=1}^{4}\mathbf{B}_{0}^{\left( \mathbf{W}_{i}\right)
}\left( t_{1},t_{2}\right) ^{2}dt_{1}dt_{2}
\end{eqnarray*}
\end{enumerate}
\end{theorem}

\bigskip

As anticipated, our proof of the above results is inspired by a proof of (%
\ref{WI}) recently given by J.-R. Pycke in \cite{Pycke}, where the author
uses a decomposition of the elements of $L^{2}\left( \left[ 0,1\right]
,dx\right) =L^{2}\left( \left[ 0,1\right] \right) $ into the orthogonal sum
of a symmetric and an antisymmetric function, around the pole $x=1/2$.
Before proving Theorem 3, we shall discuss in some detail the content of
\cite{Pycke}.

\bigskip

To this end, define for any $f\in L^{2}\left( \left[ 0,1\right] \right) $
the two operators%
\begin{equation*}
Sf\left( x\right) =\frac{1}{2}\left( f\left( x\right) +f\left( 1-x\right)
\right) \text{ \ \ and \ \ }Af\left( x\right) =\frac{1}{2}\left( f\left(
x\right) -f\left( 1-x\right) \right) \text{, \ \ }x\in \left[ 0,1\right]
\text{,}
\end{equation*}%
and observe that $f\left( x\right) =\left( A+S\right) f\left( x\right) $, $%
Sf\left( x\right) =Sf\left( 1-x\right) $ and $Af\left( x\right) =-Af\left(
1-x\right) .$ Moreover, for any $f,g\in L^{2}\left( \left[ 0,1\right]
\right) $,%
\begin{equation}
\int_{0}^{1}Af\left( x\right) Sg\left( x\right) dx=0.  \label{ortho}
\end{equation}

Note also that if $f$ is constant, then $Sf=f$ and $Af=0$.

\bigskip

\textbf{Remark -- }Let $H_{s}$ be the closed subspace of $L^{2}\left( \left[
0,1\right] \right) $ generated by functions $f$ verifying $f\left( x\right)
=f\left( 1-x\right) $ for almost every $x$, and let $H_{a}$ be the subspace
generated by functions $g$ such that $g\left( x\right) =-g\left( 1-x\right) $
for almost every $x$. Then, (\ref{ortho}) implies that $H_{s}\perp H_{a}$,
where $\perp $ indicates orthogonality in $L^{2}\left( \left[ 0,1\right]
\right) $, and also $L^{2}\left( \left[ 0,1\right] \right) =H_{s}\oplus H_{a}
$. Moreover for every $f\in L^{2}\left( \left[ 0,1\right] \right) $, $Sf$
and $Af$ equal, respectively, the orthogonal projection of $f$ on $H_{s}$,
and the orthogonal projection of $f$ on $H_{a}$.

\bigskip

The next Lemma is proved in \cite{Pycke}, and is based on a simple
computation of covariances.

\begin{lemma}
Let $b$ be a standard Brownian bridge on $\left[ 0,1\right] $, from $0$ to $%
0 $. Then, the two processes
\begin{equation*}
Ab=\left\{ Ab\left( t\right) :t\in \left[ 0,\frac{1}{2}\right] \right\}
\text{ \ \ and \ \ }Sb=\left\{ Sb\left( t\right) :t\in \left[ 0,\frac{1}{2}%
\right] \right\}
\end{equation*}%
are stochastically independent, and moreover
\begin{equation}
Ab\overset{law}{=}\left\{ \frac{b\left( 2t\right) }{2}:t\in \left[ 0,\frac{1%
}{2}\right] \right\} \text{ \ and \ }Sb\overset{law}{=}\left\{ \frac{W\left(
2t\right) }{2}:t\in \left[ 0,\frac{1}{2}\right] \right\} .  \label{identity}
\end{equation}
\end{lemma}

\bigskip

Lemma 4 yields an immediate proof of Watson's duplication identity (\ref{WI}%
). As a matter of fact, one can write, due to (\ref{ortho}) and symmetry,%
\begin{equation*}
\int_{0}^{1}\left( b\left( s\right) -\int_{0}^{1}b\left( u\right) du\right)
^{2}=2\int_{0}^{\frac{1}{2}}\left[ \left( Sb\left( t\right) -2\int_{0}^{%
\frac{1}{2}}Sb\left( u\right) du\right) ^{2}+\left( Ab\left( t\right)
\right) ^{2}\right] dt,
\end{equation*}%
and then use the relations%
\begin{equation*}
2\int_{0}^{\frac{1}{2}}\left( Ab\left( t\right) \right) ^{2}dt\overset{law}{=%
}\frac{1}{2}\int_{0}^{\frac{1}{2}}b\left( 2t\right) ^{2}dt=\frac{1}{4}%
\int_{0}^{1}b\left( v\right) ^{2}dv
\end{equation*}%
where the identity in law stems from the first part of (\ref{identity}), and%
\begin{align*}
& 2\int_{0}^{\frac{1}{2}}\left( Sb\left( t\right) -2\int_{0}^{\frac{1}{2}%
}Sb\left( u\right) du\right) ^{2}dt\overset{law}{=}\frac{1}{2}\int_{0}^{%
\frac{1}{2}}\left( W\left( 2t\right) -2\int_{0}^{\frac{1}{2}}W\left(
2u\right) du\right) ^{2}dt \\
& =\frac{1}{4}\int_{0}^{1}\left( W\left( v\right) -\int_{0}^{1}W\left(
z\right) dz\right) ^{2}dv\overset{law}{=}\frac{1}{4}\int_{0}^{1}b\left(
v\right) ^{2}dv
\end{align*}%
where the first identity in law derives again from (\ref{identity}), and the
second follows from a stochastic Fubini theorem such as the one proved e.g.
in \cite{DMartinYor}.

In the next paragraph we show that the content of Lemma 4 provides some key
elements to achieve the proof of Theorem 3.

\subsection{Proof of Theorem 3}

To prove Theorem 3 we start by defining, for every function $F$ on $\left[
0,1\right] ^{2}$, the four operators
\begin{eqnarray}
S_{1}F\left( x_{1},x_{2}\right) &=&\frac{1}{2}\left[ F\left(
x_{1},x_{2}\right) +F\left( 1-x_{1},x_{2}\right) \right]  \label{OneVOP} \\
S_{2}F\left( x_{1},x_{2}\right) &=&\frac{1}{2}\left[ F\left(
x_{1},x_{2}\right) +F\left( x_{1},1-x_{2}\right) \right]  \notag \\
A_{1}F\left( x_{1},x_{2}\right) &=&\frac{1}{2}\left[ F\left(
x_{1},x_{2}\right) -F\left( 1-x_{1},x_{2}\right) \right]  \notag \\
A_{2}F\left( x_{1},x_{2}\right) &=&\frac{1}{2}\left[ F\left(
x_{1},x_{2}\right) -F\left( x_{1},1-x_{2}\right) \right] ,  \notag
\end{eqnarray}%
where $\left( x_{1},x_{2}\right) \in \left[ 0,1\right] ^{2}$, as well as
\begin{eqnarray}
T^{\left( 1\right) }F\left( x_{1},x_{2}\right) &=&S_{1}S_{2}F\left(
x_{1},x_{2}\right) =S_{2}S_{1}F\left( x_{1},x_{2}\right)  \label{operators}
\\
T^{\left( 2\right) }F\left( x_{1},x_{2}\right) &=&S_{1}A_{2}F\left(
x_{1},x_{2}\right) =A_{2}S_{1}F\left( x_{1},x_{2}\right)  \notag \\
T^{\left( 3\right) }F\left( x_{1},x_{2}\right) &=&A_{1}S_{2}F\left(
x_{1},x_{2}\right) =S_{2}A_{1}F\left( x_{1},x_{2}\right)  \notag \\
T^{\left( 4\right) }F\left( x_{1},x_{2}\right) &=&A_{1}A_{2}F\left(
x_{1},x_{2}\right) =A_{2}A_{1}F\left( x_{1},x_{2}\right) .  \notag
\end{eqnarray}

\bigskip

Note that $F=\sum_{i=1,...,4}T^{\left( i\right) }F$, and also note the
following symmetric and antisymmetric properties: for every $\left(
x_{1},x_{2}\right) \in \left[ 0,1\right] ^{2}$,%
\begin{eqnarray*}
T^{\left( 1\right) }F\left( x_{1},x_{2}\right)  &=&T^{\left( 1\right)
}F\left( 1-x_{1},x_{2}\right) =T^{\left( 1\right) }F\left(
x_{1},1-x_{2}\right)  \\
T^{\left( 2\right) }F\left( x_{1},x_{2}\right)  &=&T^{\left( 2\right)
}F\left( 1-x_{1},x_{2}\right) =-T^{\left( 2\right) }F\left(
x_{1},1-x_{2}\right)  \\
T^{\left( 3\right) }F\left( x_{1},x_{2}\right)  &=&-T^{\left( 3\right)
}F\left( 1-x_{1},x_{2}\right) =T^{\left( 3\right) }F\left(
x_{1},1-x_{2}\right)  \\
T^{\left( 4\right) }F\left( x_{1},x_{2}\right)  &=&-T^{\left( 4\right)
}F\left( 1-x_{1},x_{2}\right) =-T^{\left( 4\right) }F\left(
x_{1},1-x_{2}\right) .
\end{eqnarray*}

This implies that, if $F$ is constant, then $T^{\left( 1\right) }F=F$, and $%
T^{\left( i\right) }F=0$ for each $i=2,3,4$. By using (\ref{ortho}) we have
moreover that, for $i\neq j$ and $F,G\in L^{2}\left( \left[ 0,1\right]
^{2},dx_{1}dx_{2}\right) =L^{2}\left( \left[ 0,1\right] ^{2}\right) $,
\begin{equation*}
\int_{\left[ 0,1\right] ^{2}}T^{\left( i\right) }F\left( x_{1},x_{2}\right)
T^{\left( j\right) }G\left( x_{1},x_{2}\right) dx_{1}dx_{2}=0,
\end{equation*}%
so that%
\begin{equation}
\int_{\left[ 0,1\right] ^{2}}F\left( x_{1},x_{2}\right)
^{2}dx_{1}dx_{2}=4\sum_{i=1}^{4}\int_{\left[ 0,\frac{1}{2}\right]
^{2}}T^{\left( i\right) }F\left( x_{1},x_{2}\right) ^{2}dx_{1}dx_{2}.
\label{norm}
\end{equation}

\bigskip

\textbf{Remark -- }Let us introduce four closed subspaces of $L^{2}\left( %
\left[ 0,1\right] ^{2}\right) $: (i) $H^{\left( 1\right) }$ is the space
generated by functions that are symmetric around the two axes $x_{1}=1/2$
and $x_{2}=1/2$; (ii) $H^{\left( 2\right) }$ is the space generated by
functions that are symmetric around the axis $x_{1}=1/2$ and antisymmetric
around $x_{2}=1/2$; (iii) $H^{\left( 3\right) }$ is the space generated by
functions $F$ that are antisymmetric around $x_{1}=1/2$ and symmetric around
$x_{2}=1/2$; (iv) $H^{\left( 4\right) }$ is the space generated by functions
$F$ that are antisymmetric around the two axes $x_{1}=1/2$ and $x_{2}=1/2$.
Then, the above relations imply that such spaces are mutually orthogonal in $%
L^{2}\left( \left[ 0,1\right] ^{2}\right) $, and $L^{2}\left( \left[ 0,1%
\right] ^{2}\right) =\oplus _{i}H^{\left( i\right) }.$ Moreover, for $%
i=1,...,4$, $T^{\left( i\right) }$, as defined in (\ref{operators}),
coincides with the orthogonal projection operator on $H^{\left( i\right) }$.
To conclude, observe that, by using standard tensor product notation%
\begin{eqnarray*}
H^{\left( 1\right) } &=&H_{s}\otimes H_{s}\text{ \ ; \ }H^{\left( 2\right)
}=H_{s}\otimes H_{a} \\
H^{\left( 3\right) } &=&H_{a}\otimes H_{s}\text{ \ ; \ }H^{\left( 4\right)
}=H_{a}\otimes H_{a},
\end{eqnarray*}%
so that $L^{2}\left( \left[ 0,1\right] ^{2}\right) =\left( H_{s}\oplus
H_{a}\right) \otimes \left( H_{s}\oplus H_{a}\right) $, where the spaces $%
H_{s}$, $H_{a}\subset L^{2}\left( \left[ 0,1\right] \right) $ have been
defined in the previous paragraph.

\subsection{Proof of part 1}

An easy calculation of covariances, based on the product formula (\ref{prod1}%
) and Lemma 4, implies that the two bivariate processes
\begin{equation*}
\left\{ A_{1}\mathbf{B}_{0}^{\left( \mathbf{W}\right) }\left(
t_{1},t_{2}\right) :\left( t_{1},t_{2}\right) \in \left[ 0,1/2\right] \times %
\left[ 0,1\right] \right\} \text{ \ and \ }\left\{ S_{1}\mathbf{B}%
_{0}^{\left( \mathbf{W}\right) }\left( t_{1},t_{2}\right) :\left(
t_{1},t_{2}\right) \in \left[ 0,1/2\right] \times \left[ 0,1\right] \right\}
\end{equation*}%
are stochastically independent, and an analogous conclusion holds for the
two processes
\begin{equation*}
\left\{ A_{2}\mathbf{B}_{0}^{\left( \mathbf{W}\right) }\left(
t_{1},t_{2}\right) :\left( t_{1},t_{2}\right) \in \left[ 0,1\right] \times %
\left[ 0,1/2\right] \right\} \text{ \ and \ }\left\{ S_{2}\mathbf{B}%
_{0}^{\left( \mathbf{W}\right) }\left( t_{1},t_{2}\right) :\left(
t_{1},t_{2}\right) \in \left[ 0,1\right] \times \left[ 0,1/2\right] \right\}
.
\end{equation*}

This entails immediately that the four (jointly) Gaussian processes
\begin{equation*}
\left\{ T^{\left( i\right) }\mathbf{B}_{0}^{\left( \mathbf{W}\right) }\left(
t_{1},t_{2}\right) :\left( t_{1},t_{2}\right) \in \left[ 0,1/2\right]
^{2}\right\} \text{, \ }i=1,...,4\text{,}
\end{equation*}%
are mutually independent. Now, by applying (\ref{norm}) to the random
continuous function
\begin{equation*}
\left( t_{1},t_{2}\right) \mapsto \mathbf{B}_{0}^{\left( \mathbf{W}\right)
}\left( t_{1},t_{2}\right) -\int_{\left[ 0,1\right] ^{2}}\mathbf{B}%
_{0}^{\left( \mathbf{W}\right) }\left( u_{1},u_{2}\right) du_{1}du_{2}
\end{equation*}%
we obtain, thanks to symmetry,%
\begin{eqnarray*}
&&\int_{\left[ 0,1\right] ^{2}}\left[ \mathbf{B}_{0}^{\left( \mathbf{W}%
\right) }\left( t_{1},t_{2}\right) -\int_{\left[ 0,1\right] ^{2}}du_{1}du_{2}%
\mathbf{B}_{0}^{\left( \mathbf{W}\right) }\left( u_{1},u_{2}\right) \right]
^{2}dt_{1}dt_{2} \\
&=&4\int_{\left[ 0,\frac{1}{2}\right] ^{2}}\left[ T^{\left( 1\right) }%
\mathbf{B}_{0}^{\left( \mathbf{W}\right) }\left( t_{1},t_{2}\right) -4\int_{%
\left[ 0,\frac{1}{2}\right] ^{2}}T^{\left( 1\right) }\mathbf{B}_{0}^{\left(
\mathbf{W}\right) }\left( u_{1},u_{2}\right) du_{1}du_{2}\right]
^{2}dt_{1}dt_{2} \\
&&+4\sum_{i=2}^{4}\int_{\left[ 0,\frac{1}{2}\right] ^{2}}T^{\left( i\right) }%
\mathbf{B}_{0}^{\left( \mathbf{W}\right) }\left( t_{1},t_{2}\right)
^{2}dt_{1}dt_{2}.
\end{eqnarray*}

Since for any Brownian sheet $\mathbf{W}$%
\begin{equation*}
\mathbf{K}^{\left( \mathbf{W},1\right) }\left( t_{2},t_{1}\right) \overset{%
law}{=}\mathbf{K}^{\left( \mathbf{W},2\right) }\left( t_{1},t_{2}\right)
\end{equation*}%
where the identity holds for the two processes as a whole, the proof of
Theorem 3 is achieved once the following three identities in law are shown,%
\begin{align}
& 4\int_{\left[ 0,\frac{1}{2}\right] ^{2}}\left[ T^{\left( 1\right) }\mathbf{%
B}_{0}^{\left( \mathbf{W}\right) }\left( t_{1},t_{2}\right) -4\int_{\left[ 0,%
\frac{1}{2}\right] ^{2}}T^{\left( 1\right) }\mathbf{B}_{0}^{\left( \mathbf{W}%
\right) }\left( u_{1},u_{2}\right) du_{1}du_{2}\right] ^{2}dt_{1}dt_{2}
\label{F1} \\
& \overset{law}{=}\frac{1}{16}\int_{\left[ 0,1\right] ^{2}}\mathbf{B}%
^{\left( \mathbf{W}\right) }\left( t_{1},t_{2}\right) ^{2}dt_{1}dt_{2}
\notag
\end{align}%
\begin{align}
& 4\int_{\left[ 0,\frac{1}{2}\right] ^{2}}T^{\left( 2\right) }\mathbf{B}%
_{0}^{\left( \mathbf{W}\right) }\left( t_{1},t_{2}\right) ^{2}dt_{1}dt_{2}%
\overset{law}{=}\frac{1}{16}\int_{\left[ 0,1\right] ^{2}}\mathbf{K}^{\left(
\mathbf{W},1\right) }\left( t_{1},t_{2}\right) ^{2}dt_{1}dt_{2}  \label{F2}
\\
& \overset{law}{=}4\int_{\left[ 0,\frac{1}{2}\right] ^{2}}T^{\left( 3\right)
}\mathbf{B}_{0}^{\left( \mathbf{W}\right) }\left( t_{1},t_{2}\right)
^{2}dt_{1}dt_{2}  \notag
\end{align}%
\begin{equation}
4\int_{\left[ 0,\frac{1}{2}\right] ^{2}}T^{\left( 4\right) }\mathbf{B}%
_{0}^{\left( \mathbf{W}\right) }\left( t_{1},t_{2}\right) ^{2}dt_{1}dt_{2}%
\overset{law}{=}\frac{1}{16}\int_{\left[ 0,1\right] ^{2}}\mathbf{B}%
_{0}^{\left( \mathbf{W}\right) }\left( t_{1},t_{2}\right) ^{2}dt_{1}dt_{2}
\label{F3}
\end{equation}

To prove (\ref{F1}), just observe that Lemma 4 and (\ref{prod3}) entail
\begin{equation}
\left\{ S_{1}\mathbf{B}_{0}^{\left( \mathbf{W}\right) }\left(
t_{1},t_{2}\right) :\left( t_{1},t_{2}\right) \in \left[ 0,1/2\right] \times %
\left[ 0,1\right] \right\} \overset{law}{=}\left\{ 2^{-\frac{1}{2}}\mathbf{K}%
^{\left( \mathbf{W},2\right) }\left( t_{1},t_{2}\right) :\left(
t_{1},t_{2}\right) \in \left[ 0,1/2\right] \times \left[ 0,1\right] \right\}
\label{S-identity}
\end{equation}%
and therefore%
\begin{equation}
\left\{ T^{\left( 1\right) }\mathbf{B}_{0}^{\left( \mathbf{W}\right) }\left(
t_{1},t_{2}\right) :\left( t_{1},t_{2}\right) \in \left[ 0,1/2\right]
^{2}\right\} \overset{law}{=}\left\{ 2^{-1}\mathbf{W}\left(
t_{1},t_{2}\right) :\left( t_{1},t_{2}\right) \in \left[ 0,1/2\right]
^{2}\right\}  \label{B-W}
\end{equation}%
so that%
\begin{align*}
& 4\int_{\left[ 0,\frac{1}{2}\right] ^{2}}\left[ T^{\left( 1\right) }\mathbf{%
B}_{0}^{\left( \mathbf{W}\right) }\left( t_{1},t_{2}\right) -4\int_{\left[ 0,%
\frac{1}{2}\right] ^{2}}T^{\left( 1\right) }\mathbf{B}_{0}^{\left( \mathbf{W}%
\right) }\left( u_{1},u_{2}\right) du_{1}du_{2}\right] ^{2}dt_{1}dt_{2} \\
& \overset{law}{=}\int_{\left[ 0,\frac{1}{2}\right] ^{2}}\left[ \mathbf{W}%
\left( t_{1},t_{2}\right) -4\int_{\left[ 0,\frac{1}{2}\right] ^{2}}\mathbf{W}%
\left( u_{1},u_{2}\right) du_{1}du_{2}\right] ^{2}dt_{1}dt_{2} \\
& \overset{law}{=}\frac{1}{4}\int_{\left[ 0,\frac{1}{2}\right] ^{2}}\left[
\mathbf{W}\left( 2t_{1},2t_{2}\right) -4\int_{\left[ 0,\frac{1}{2}\right]
^{2}}\mathbf{W}\left( 2u_{1},2u_{2}\right) du_{1}du_{2}\right]
^{2}dt_{1}dt_{2} \\
& =\frac{1}{16}\int_{\left[ 0,1\right] ^{2}}\left[ \mathbf{W}\left(
s_{1},s_{2}\right) -\int_{\left[ 0,1\right] ^{2}}\mathbf{W}\left(
v_{1},v_{2}\right) dv_{1}dv_{2}\right] ^{2}ds_{1}ds_{2} \\
& \overset{law}{=}\frac{1}{16}\int_{\left[ 0,1\right] ^{2}}\mathbf{B}%
^{\left( \mathbf{W}\right) }\left( s_{1},s_{2}\right) ^{2}ds_{1}ds_{2}
\end{align*}%
where the last equality is a consequence of a stochastic Fubini theorem, and
namely of relation (\ref{Fub1}) in the statement of Corollary 2.

To prove (\ref{F2}), we use (\ref{S-identity}), (\ref{prod3}) and Lemma 4 to
obtain that
\begin{equation*}
\left\{ T^{\left( 2\right) }\mathbf{B}_{0}^{\left( \mathbf{W}\right) }\left(
t_{1},t_{2}\right) :\left( t_{1},t_{2}\right) \in \left[ 0,1/2\right]
^{2}\right\} \overset{law}{=}\left\{ 2^{-\frac{3}{2}}\mathbf{K}^{\left(
\mathbf{W},2\right) }\left( t_{1},2t_{2}\right) :\left( t_{1},t_{2}\right)
\in \left[ 0,1/2\right] ^{2}\right\}
\end{equation*}%
and eventually%
\begin{gather*}
4\int_{\left[ 0,\frac{1}{2}\right] ^{2}}T^{\left( 2\right) }\mathbf{B}%
_{0}^{\left( \mathbf{W}\right) }\left( t_{1},t_{2}\right) ^{2}dt_{1}dt_{2}%
\overset{law}{=}\frac{1}{2}\int_{\left[ 0,\frac{1}{2}\right] ^{2}}\mathbf{K}%
^{\left( \mathbf{W},2\right) }\left( t_{1},2t_{2}\right) ^{2}dt_{1}dt_{2} \\
\overset{law}{=}\frac{1}{4}\int_{\left[ 0,\frac{1}{2}\right] ^{2}}\mathbf{K}%
^{\left( \mathbf{W},2\right) }\left( 2t_{1},2t_{2}\right) ^{2}dt_{1}dt_{2}=%
\frac{1}{16}\int_{\left[ 0,1\right] ^{2}}\mathbf{K}^{\left( \mathbf{W}%
,2\right) }\left( u_{1},u_{2}\right) ^{2}du_{1}du_{2}.
\end{gather*}

The case of $T^{\left( 3\right) }$ can be treated analogously by using (\ref%
{prod2}). To conclude, we note that
\begin{equation*}
\left\{ A_{1}\mathbf{B}_{0}^{\left( \mathbf{W}\right) }\left(
t_{1},t_{2}\right) :\left( t_{1},t_{2}\right) \in \left[ 0,1/2\right] \times %
\left[ 0,1\right] \right\} \overset{law}{=}\left\{ 2^{-1}\mathbf{B}%
_{0}^{\left( \mathbf{W}\right) }\left( 2t_{1},t_{2}\right) :\left(
t_{1},t_{2}\right) \in \left[ 0,1/2\right] \times \left[ 0,1\right] \right\}
\end{equation*}%
and therefore%
\begin{equation*}
\left\{ T^{\left( 4\right) }\mathbf{B}_{0}^{\left( \mathbf{W}\right) }\left(
t_{1},t_{2}\right) :\left( t_{1},t_{2}\right) \in \left[ 0,1/2\right]
^{2}\right\} \overset{law}{=}\left\{ 2^{-2}\mathbf{B}_{0}^{\left( \mathbf{W}%
\right) }\left( 2t_{1},2t_{2}\right) :\left( t_{1},t_{2}\right) \in \left[
0,1/2\right] ^{2}\right\} ,
\end{equation*}%
so that%
\begin{equation*}
4\int_{\left[ 0,\frac{1}{2}\right] ^{2}}T^{\left( 4\right) }\mathbf{B}%
_{0}^{\left( \mathbf{W}\right) }\left( t_{1},t_{2}\right) ^{2}dt_{1}dt_{2}%
\overset{law}{=}\frac{1}{4}\int_{\left[ 0,\frac{1}{2}\right] ^{2}}\mathbf{B}%
_{0}^{\left( \mathbf{W}\right) }\left( 2t_{1},2t_{2}\right) ^{2}dt_{1}dt_{2}=%
\frac{1}{16}\int_{\left[ 0,1\right] ^{2}}\mathbf{B}_{0}^{\left( \mathbf{W}%
\right) }\left( u_{1},u_{2}\right) ^{2}du_{1}du_{2}
\end{equation*}

\subsection{Proof of part 2}

We write%
\begin{equation*}
\mathbf{B}_{0}^{\left( \mathbf{W}\right) }\left( t_{1},t_{2}\right)
-\int_{0}^{1}\mathbf{B}_{0}^{\left( \mathbf{W}\right) }\left(
t_{1},u_{2}\right) du_{2}=S_{2}\mathbf{B}_{0}^{\left( \mathbf{W}\right)
}\left( t_{1},t_{2}\right) -\int_{0}^{1}S_{2}\mathbf{B}_{0}^{\left( \mathbf{W%
}\right) }\left( t_{1},u_{2}\right) du_{2}+A_{2}\mathbf{B}_{0}^{\left(
\mathbf{W}\right) }\left( t_{1},t_{2}\right) ,
\end{equation*}%
where the operators $S_{2}$ and $A_{2}$ are defined in (\ref{OneVOP}). Since
$S_{2}=T^{\left( 1\right) }+T^{\left( 3\right) }$ and $A_{2}=T^{\left(
2\right) }+T^{\left( 4\right) }$, we can use orthogonality and symmetry to
obtain%
\begin{eqnarray*}
&&\int_{\left[ 0,1\right] ^{2}}\left[ \mathbf{B}_{0}^{\left( \mathbf{W}%
\right) }\left( t_{1},t_{2}\right) -\int_{0}^{1}\mathbf{B}_{0}^{\left(
\mathbf{W}\right) }\left( t_{1},u_{2}\right) du_{2}\right] ^{2}dt_{1}dt_{2}
\\
&=&\int_{\left[ 0,1\right] ^{2}}\left[ S_{2}\mathbf{B}_{0}^{\left( \mathbf{W}%
\right) }\left( t_{1},t_{2}\right) -\int_{0}^{1}S_{2}\mathbf{B}_{0}^{\left(
\mathbf{W}\right) }\left( t_{1},u_{2}\right) du_{2}\right] ^{2}dt_{1}dt_{2}
\\
&&+\int_{\left[ 0,1\right] ^{2}}A_{2}\mathbf{B}_{0}^{\left( \mathbf{W}%
\right) }\left( t_{1},t_{2}\right) ^{2}dt_{1}dt_{2} \\
&=&2\int_{\left[ 0,1\right] \times \left[ 0,1/2\right] }\left[ S_{2}\mathbf{B%
}_{0}^{\left( \mathbf{W}\right) }\left( t_{1},t_{2}\right) -2\int_{0}^{\frac{%
1}{2}}S_{2}\mathbf{B}_{0}^{\left( \mathbf{W}\right) }\left(
t_{1},u_{2}\right) du_{2}\right] ^{2}dt_{1}dt_{2} \\
&&+2\int_{\left[ 0,1\right] \times \left[ 0,1/2\right] }A_{2}\mathbf{B}%
_{0}^{\left( \mathbf{W}\right) }\left( t_{1},t_{2}\right) ^{2}dt_{1}dt_{2}.
\end{eqnarray*}

We already know that the restrictions to $\left[ 0,1\right] \times \left[
0,1/2\right] $ of the two processes $S_{2}\mathbf{B}_{0}^{\left( \mathbf{W}%
\right) }$ and $A_{2}\mathbf{B}_{0}^{\left( \mathbf{W}\right) }$ are
stochastically independent. Moreover Lemma 4 and (\ref{prod2}) imply the two
relations
\begin{eqnarray}
&&\left\{ S_{2}\mathbf{B}_{0}^{\left( \mathbf{W}\right) }\left(
t_{1},t_{2}\right) :\left( t_{1},t_{2}\right) \in \left[ 0,1\right] \times %
\left[ 0,1/2\right] \right\}  \label{idProc} \\
&&\overset{law}{=}\left\{ 2^{-\frac{1}{2}}\mathbf{K}^{\left( \mathbf{W}%
,1\right) }\left( t_{1},t_{2}\right) :\left( t_{1},t_{2}\right) \in \left[
0,1\right] \times \left[ 0,1/2\right] \right\}  \notag \\
&&\left\{ A_{2}\mathbf{B}_{0}^{\left( \mathbf{W}\right) }\left(
t_{1},t_{2}\right) :\left( t_{1},t_{2}\right) \in \left[ 0,1\right] \times %
\left[ 0,1/2\right] \right\}  \notag \\
&&\overset{law}{=}\left\{ 2^{-1}\mathbf{B}_{0}^{\left( \mathbf{W}\right)
}\left( t_{1},2t_{2}\right) :\left( t_{1},t_{2}\right) \in \left[ 0,1\right]
\times \left[ 0,1/2\right] \right\} .  \notag
\end{eqnarray}

As a consequence, we obtain
\begin{eqnarray*}
&&2\int_{\left[ 0,1\right] \times \left[ 0,1/2\right] }A_{2}\mathbf{B}%
_{0}^{\left( \mathbf{W}\right) }\left( t_{1},t_{2}\right) ^{2}dt_{1}dt_{2}%
\overset{law}{=}\frac{1}{2}\int_{\left[ 0,1\right] \times \left[ 0,1/2\right]
}\mathbf{B}_{0}^{\left( \mathbf{W}\right) }\left( t_{1},2t_{2}\right)
^{2}dt_{1}dt_{2} \\
&=&\frac{1}{4}\int_{\left[ 0,1\right] ^{2}}\mathbf{B}_{0}^{\left( \mathbf{W}%
\right) }\left( t_{1},t_{2}\right) ^{2}dt_{1}dt_{2}.
\end{eqnarray*}

To conclude the proof, use the first part of (\ref{idProc}) and scaling to
obtain%
\begin{align*}
& 2\int_{\left[ 0,1\right] \times \left[ 0,1/2\right] }\left[ S_{2}\mathbf{B}%
_{0}^{\left( \mathbf{W}\right) }\left( t_{1},t_{2}\right) -2\int_{0}^{\frac{1%
}{2}}S_{2}\mathbf{B}_{0}^{\left( \mathbf{W}\right) }\left(
t_{1},u_{2}\right) du_{2}\right] ^{2}dt_{1}dt_{2} \\
& \overset{law}{=}\frac{1}{2}\int_{\left[ 0,1\right] \times \left[ 0,1/2%
\right] }\left[ \mathbf{K}^{\left( \mathbf{W},1\right) }\left(
t_{1},2t_{2}\right) -2\int_{0}^{\frac{1}{2}}\mathbf{K}^{\left( \mathbf{W}%
,1\right) }\left( t_{1},2u_{2}\right) du_{2}\right] ^{2}dt_{1}dt_{2} \\
& =\frac{1}{4}\int_{\left[ 0,1\right] ^{2}}\left[ \mathbf{K}^{\left( \mathbf{%
W},1\right) }\left( t_{1},t_{2}\right) -\int_{0}^{1}\mathbf{K}^{\left(
\mathbf{W},1\right) }\left( t_{1},u_{2}\right) du_{2}\right]
^{2}dt_{1}dt_{2}.
\end{align*}

Now define $\left\{ \lambda _{i},f_{i}:i\geq 1\right\} $ and $\left\{ \gamma
_{i},g_{i}:i\geq 1\right\} $ to be the sequences of eigenvalues and
eigenfunctions of the Hilbert-Schmidt operators associated to the covariance
function, respectively of $t\mapsto b\left( t\right) $, and of
\begin{equation*}
t\mapsto Z\left( t\right) :=W\left( t\right) -\int_{0}^{1}W\left( z\right)
dz.
\end{equation*}

It is well known (see e.g. \cite{Kac1}) that there exist two sequences $%
\left\{ \xi _{i}:i\geq 1\right\} $ and $\left\{ \zeta _{i}:i\geq 1\right\} $
of i.i.d. standard Gaussian random variables such that the Karhunen-Lo\`{e}%
ve expansions of $b$ and $Z$ are respectively given by%
\begin{equation*}
b\left( t\right) =\sum_{i\geq 1}\xi _{i}\sqrt{\lambda _{i}}f_{i}\left(
t\right) \text{ \ \ and \ \ }Z\left( t\right) =\sum_{i\geq 1}\zeta _{i}\sqrt{%
\gamma _{i}}g_{i}\left( t\right) ,
\end{equation*}%
and moreover (see \cite{DMartinYor}) $\gamma _{i}=\lambda _{i}$ for every $%
i\geq 1$. Since (\ref{prod2}) implies that, for every $\left(
t_{1},t_{2}\right) ,\left( s_{1},s_{2}\right) \in \left[ 0,1\right] ^{2}$,%
\begin{eqnarray*}
&&\mathbb{E}\left[ \left( \mathbf{K}^{\left( \mathbf{W},1\right) }\left(
t_{1},t_{2}\right) -\int_{0}^{1}\mathbf{K}^{\left( \mathbf{W},1\right)
}\left( t_{1},u_{2}\right) du_{2}\right) \left( \mathbf{K}^{\left( \mathbf{W}%
,1\right) }\left( s_{1},s_{2}\right) -\int_{0}^{1}\mathbf{K}^{\left( \mathbf{%
W},1\right) }\left( s_{1},u_{2}\right) du_{2}\right) \right] \\
&=&\mathbb{E}\left[ b\left( t_{1}\right) b\left( s_{1}\right) \right] \times
\mathbb{E}\left[ \left( W\left( t_{2}\right) -\int_{0}^{1}W\left( z\right)
dz\right) \left( W\left( s_{2}\right) -\int_{0}^{1}W\left( z\right)
dz\right) \right]
\end{eqnarray*}%
we conclude immediately (by using, for instance, \cite[Lemma 4.1]{DPY}) that
the Karhunen-Lo\`{e}ve expansion of the bivariate Gaussian process
\begin{equation*}
\mathbf{Z}\left( s,t\right) =\mathbf{K}^{\left( \mathbf{W},1\right) }\left(
s,t\right) -\int_{0}^{1}\mathbf{K}^{\left( \mathbf{W},1\right) }\left(
s,u\right) du
\end{equation*}%
is given by
\begin{equation*}
\mathbf{Z}\left( s,t\right) =\sum_{i,j\geq 1}\sqrt{\lambda _{i}\lambda _{j}}%
\theta _{ij}f_{i}\left( s\right) g_{i}\left( t\right) \text{,}
\end{equation*}%
where $\left\{ \theta _{ij}:i,j\geq 1\right\} $ is an array of i.i.d.
standard Gaussian random variables. This last relation entails that%
\begin{eqnarray}
&&\frac{1}{4}\int_{\left[ 0,1\right] ^{2}}\left[ \mathbf{K}^{\left( \mathbf{W%
},1\right) }\left( t_{1},t_{2}\right) -\int_{0}^{1}\mathbf{K}^{\left(
\mathbf{W},1\right) }\left( t_{1},u_{2}\right) du_{2}\right] ^{2}dt_{1}dt_{2}
\label{KLidentity} \\
&=&\frac{1}{4}\sum_{i\geq 1}\lambda _{i}\sum_{j\geq 1}\lambda _{j}\theta
_{ij}^{2}\overset{law}{=}\frac{1}{4}\int_{\left[ 0,1\right] ^{2}}\mathbf{B}%
_{0}^{\left( \mathbf{W}\right) }\left( t_{1},t_{2}\right) ^{2}dt_{1}dt_{2}%
\text{.}  \notag
\end{eqnarray}

To justify the last equality in law, just observe that, thanks again to \cite%
[Lemma 4.1]{DPY} and formula (\ref{prod1}), the Karhunen-Lo\`{e}ve expansion
of $\mathbf{B}_{0}$ is given by
\begin{equation*}
\sum_{i,j\geq 1}\sqrt{\lambda _{i}\lambda _{j}}\eta _{ij}g_{i}\left(
s\right) g_{i}\left( t\right)
\end{equation*}%
where $\left\{ \eta _{ij}:i,j\geq 1\right\} $ is an array of i.i.d. standard
Gaussian random variables (the reader is referred to \cite{DPY} for a
detailed discussion of Karhunen-Lo\`{e}ve expansions for bivariate Gaussian
processes).

\subsection{Proof of part 3}

We first observe that
\begin{align*}
0& =T^{\left( 2\right) }\int_{0}^{1}\mathbf{B}_{0}^{\left( \mathbf{W}\right)
}\left( t_{1},u_{2}\right) du_{2}=T^{\left( 4\right) }\int_{0}^{1}\mathbf{B}%
_{0}^{\left( \mathbf{W}\right) }\left( t_{1},u_{2}\right) du_{2} \\
& =T^{\left( 3\right) }\int_{0}^{1}\mathbf{B}_{0}^{\left( \mathbf{W}\right)
}\left( u_{1},t_{2}\right) du_{1}=T^{\left( 4\right) }\int_{0}^{1}\mathbf{B}%
_{0}^{\left( \mathbf{W}\right) }\left( u_{1},t_{2}\right) du_{1} \\
& =T^{\left( i\right) }\int_{0}^{1}\mathbf{B}_{0}^{\left( \mathbf{W}\right)
}\left( u_{1},u_{2}\right) du_{1}du_{2},\text{ \ \ }i=2,3,4,
\end{align*}%
and
\begin{eqnarray*}
T^{\left( i\right) }\int_{0}^{1}\mathbf{B}_{0}^{\left( \mathbf{W}\right)
}\left( u_{1},t_{2}\right) du_{1} &=&\int_{0}^{1}T^{\left( i\right) }\mathbf{%
B}_{0}^{\left( \mathbf{W}\right) }\left( u_{1},t_{2}\right) du_{1}\text{, \
\ }i=1,2 \\
T^{\left( j\right) }\int_{0}^{1}\mathbf{B}_{0}^{\left( \mathbf{W}\right)
}\left( t_{1},u_{2}\right) du_{2} &=&\int_{0}^{1}T^{\left( j\right) }\mathbf{%
B}_{0}^{\left( \mathbf{W}\right) }\left( t_{1},u_{2}\right) du_{2}\text{, \
\ }j=1,3.
\end{eqnarray*}

As a consequence, by orthogonality and symmetry,%
\begin{align*}
& \int_{\left[ 0,1\right] ^{2}}\left[ \mathbf{B}_{0}^{\left( \mathbf{W}%
\right) }\left( t_{1},t_{2}\right) -\int_{0}^{1}\mathbf{B}_{0}^{\left(
\mathbf{W}\right) }\left( t_{1},u_{2}\right) du_{2}\right. \\
& \left. -\int_{0}^{1}\mathbf{B}_{0}^{\left( \mathbf{W}\right) }\left(
u_{1},t_{2}\right) du_{1}+\int_{\left[ 0,\frac{1}{2}\right] ^{2}}\mathbf{B}%
_{0}^{\left( \mathbf{W}\right) }\left( u_{1},u_{2}\right) du_{1}du_{2}\right]
^{2}dt_{1}dt_{2} \\
& =4\int_{\left[ 0,\frac{1}{2}\right] ^{2}}\left[ T^{\left( 1\right) }%
\mathbf{B}_{0}^{\left( \mathbf{W}\right) }\left( t_{1},t_{2}\right)
-2\int_{0}^{\frac{1}{2}}T^{\left( 1\right) }\mathbf{B}_{0}^{\left( \mathbf{W}%
\right) }\left( t_{1},u_{2}\right) du_{2}\right. \\
& \left. -2\int_{0}^{\frac{1}{2}}T^{\left( 1\right) }\mathbf{B}_{0}^{\left(
\mathbf{W}\right) }\left( u_{1},t_{2}\right) du_{1}+4\int_{\left[ 0,\frac{1}{%
2}\right] ^{2}}T^{\left( 1\right) }\mathbf{B}_{0}^{\left( \mathbf{W}\right)
}\left( u_{1},u_{2}\right) du_{1}du_{2}\right] ^{2}dt_{1}dt_{2} \\
& +4\int_{\left[ 0,\frac{1}{2}\right] ^{2}}\left[ T^{\left( 2\right) }%
\mathbf{B}_{0}^{\left( \mathbf{W}\right) }\left( t_{1},t_{2}\right)
-2\int_{0}^{\frac{1}{2}}T^{\left( 2\right) }\mathbf{B}_{0}^{\left( \mathbf{W}%
\right) }\left( u_{1},t_{2}\right) du_{1}\right] ^{2}dt_{1}dt_{2} \\
& +4\int_{\left[ 0,\frac{1}{2}\right] ^{2}}\left[ T^{\left( 3\right) }%
\mathbf{B}_{0}^{\left( \mathbf{W}\right) }\left( t_{1},t_{2}\right)
-2\int_{0}^{\frac{1}{2}}T^{\left( 3\right) }\mathbf{B}_{0}^{\left( \mathbf{W}%
\right) }\left( t_{1},u_{2}\right) du_{2}\right] ^{2}dt_{1}dt_{2} \\
& +4\int_{\left[ 0,\frac{1}{2}\right] ^{2}}\left[ T^{\left( 4\right) }%
\mathbf{B}_{0}^{\left( \mathbf{W}\right) }\left( t_{1},t_{2}\right) \right]
^{2}dt_{1}dt_{2} \\
& \overset{def}{=}Q_{1}+Q_{2}+Q_{3}+Q_{4}.
\end{align*}

Since we know, thanks to the previous discussion, that the $Q_{i}$'s are
mutually independent, it is now sufficient to show that, for $i=1,...,4$,
\begin{equation}
Q_{i}\overset{law}{=}\frac{1}{16}\int_{\left[ 0,1\right] ^{2}}\mathbf{B}%
_{0}^{\left( \mathbf{W}_{i}\right) }\left( t_{1},t_{2}\right)
^{2}dt_{1}dt_{2}.  \label{Eqdistr}
\end{equation}

We start with $Q_{2}$ (by symmetry, the case of $Q_{3}$ is handled
analogously), and recall that we have already proved that%
\begin{eqnarray*}
&&Q_{i}\overset{law}{=}\frac{1}{4}\int_{\left[ 0,\frac{1}{2}\right] ^{2}}%
\left[ \mathbf{K}^{\left( \mathbf{W},2\right) }\left( 2t_{1},2t_{2}\right)
-2\int_{0}^{\frac{1}{2}}\mathbf{K}^{\left( \mathbf{W},2\right) }\left(
2u_{1},2t_{2}\right) du_{1}\right] ^{2}dt_{1}dt_{2} \\
&=&\frac{1}{16}\int_{\left[ 0,1\right] ^{2}}\left[ \mathbf{K}^{\left(
\mathbf{W},2\right) }\left( v_{1},v_{2}\right) -\int_{0}^{1}\mathbf{K}%
^{\left( \mathbf{W},2\right) }\left( z,v_{2}\right) dz\right]
^{2}dv_{1}dv_{2}
\end{eqnarray*}%
so that (\ref{Eqdistr}) in the case $i=2,3$ derives immediately from (\ref%
{KLidentity}). Since we have proved (\ref{Eqdistr}) for $i=4$ (to obtain
part 1 of Theorem 3) we are now left with the case $i=1$.

To see that (\ref{Eqdistr}) holds also in this case, use (\ref{B-W}) to
write, after a standard change of variables,%
\begin{align*}
& Q_{1}\overset{law}{=}\frac{1}{16}\int_{\left[ 0,1\right] ^{2}}\left[
\mathbf{W}\left( t_{1},t_{2}\right) -\int_{0}^{1}\mathbf{W}\left(
t_{1},u_{2}\right) du_{2}\right.  \\
& \left. -\int_{0}^{1}\mathbf{W}\left( u_{1},t_{2}\right) du_{1}+\int_{\left[
0,1\right] ^{2}}\mathbf{W}\left( u_{1},u_{2}\right) du_{1}du_{2}\right]
^{2}dt_{1}dt_{2}
\end{align*}%
and then apply relation (\ref{Fub2}) in Corollary 2.

\bigskip

\textbf{Remark -- }Note that the techniques used for the proof of Therorem 3
could be also applied to the case of general $n$-variate Gaussian processes,
for $n>2$.

\section{Application: Fourier transforms of double Wiener integrals with
respect to conditioned Gaussian processes}

Let the above notation prevail, and let $\mathbf{W}_{1}$ and $\mathbf{W}_{2}$
be two independent Brownian sheets. In this section, we are interested in
finding the explicit Fourier transform of the three double Wiener integrals
\begin{eqnarray*}
\mathbf{I} &=&\int_{\left[ 0,1\right] ^{2}}\mathbf{B}_{0}^{\left( \mathbf{W}%
_{1}\right) }\left( t_{1},t_{2}\right) \mathbf{B}^{\left( \mathbf{W}%
_{2}\right) }\left( dt_{1},dt_{2}\right) \\
&=&\int_{\left[ 0,1\right] ^{2}}\left[ \mathbf{B}_{0}^{\left( \mathbf{W}%
_{1}\right) }\left( t_{1},t_{2}\right) -\int_{\left[ 0,1\right] ^{2}}\mathbf{%
B}_{0}^{\left( \mathbf{W}_{1}\right) }\left( u_{1},u_{2}\right) du_{1}du_{2}%
\right] \mathbf{B}^{\left( \mathbf{W}_{2}\right) }\left(
dt_{1},dt_{2}\right) ; \\
\mathbf{J} &=&\int_{\left[ 0,1\right] ^{2}}\mathbf{B}_{0}^{\left( \mathbf{W}%
_{1}\right) }\left( t_{1},t_{2}\right) \mathbf{K}^{\left( \mathbf{W}%
_{2},2\right) }\left( dt_{1},dt_{2}\right) \\
&=&\int_{\left[ 0,1\right] ^{2}}\left[ \mathbf{B}_{0}^{\left( \mathbf{W}%
_{1}\right) }\left( t_{1},t_{2}\right) -\int_{\left[ 0,1\right] ^{2}}\mathbf{%
B}_{0}^{\left( \mathbf{W}_{1}\right) }\left( t_{1},u_{2}\right) du_{2}\right]
\mathbf{K}^{\left( \mathbf{W}_{2},2\right) }\left( dt_{1},dt_{2}\right) ; \\
\mathbf{Y} &\mathbf{=}&\int_{\left[ 0,1\right] ^{2}}\mathbf{B}_{0}^{\left(
\mathbf{W}_{1}\right) }\left( t_{1},t_{2}\right) \mathbf{B}_{0}^{\left(
\mathbf{W}_{2}\right) }\left( dt_{1},dt_{2}\right) \\
&=&\int_{\left[ 0,1\right] ^{2}}\left[ \mathbf{B}_{0}^{\left( \mathbf{W}%
_{1}\right) }\left( t_{1},t_{2}\right) -\int_{0}^{1}\mathbf{B}_{0}^{\left(
\mathbf{W}_{1}\right) }\left( t_{1},u_{2}\right) du_{2}\right. \\
&&\left. -\int_{0}^{1}\mathbf{B}_{0}^{\left( \mathbf{W}_{1}\right) }\left(
u_{1},t_{2}\right) du_{1}+\int_{\left[ 0,1\right] ^{2}}\mathbf{B}%
_{0}^{\left( \mathbf{W}_{1}\right) }\left( u_{1},u_{2}\right) du_{1}du_{2}%
\right] \mathbf{B}_{0}^{\left( \mathbf{W}_{2}\right) }\left(
dt_{1},dt_{2}\right)
\end{eqnarray*}

We shall show that such computations can be achieved by means of Theorem 3.
To this end, we introduce some notation borrowed from \cite{DPY}: for every $%
a\in \mathbb{C}$,

\begin{enumerate}
\item $C\left( a\right) =\prod_{j=1}^{\infty }\cosh \left( \frac{a}{j\pi }%
\right) $;

\item $C_{\text{odd}}\left( a\right) =\prod_{j=0}^{\infty }\cosh \left[
\frac{a}{\left( 2j+1\right) \pi }\right] $;

\item $C_{\text{even}}\left( a\right) =\prod_{j=1}^{\infty }\cosh \left[
\frac{a}{2j\pi }\right] =C\left( \frac{a}{2}\right) ;$

\item $S\left( a\right) =\prod_{j=1}^{\infty }\left[ \pi j\sinh \left( \frac{%
a}{\pi j}\right) /a\right] ;$

\item $S_{\text{even}}\left( a\right) =\prod_{j=1}^{\infty }\left[ \pi
2j\sinh \left( \frac{a}{\pi 2j}\right) /a\right] =S\left( a/2\right) ;$

\item $S_{\text{odd}}\left( a\right) =\prod_{j=1}^{\infty }\left[ \pi \left(
2j-1\right) \sinh \left( \frac{a}{\pi \left( 2j-1\right) }\right) /a\right]
=C\left( a/2\right) $;

\item $\mathcal{T}\left( a\right) =\sum_{j=0}^{\infty }\left\{ \tanh \left(
\frac{2a}{\left( 2j+1\right) \pi }\right) \left[ \left( 2j+1\right) \pi %
\right] ^{-1}\right\} $.
\end{enumerate}

\bigskip

Moreover, we recall the following result

\begin{proposition}[{see \protect\cite[Proposition 4.1]{DPY}}]
For every $u\in \Re $

\begin{enumerate}
\item $\mathbb{E}\left[ \exp \left( -\frac{u^{2}}{2}\int_{\left[ 0,1\right]
^{2}}\mathbf{B}^{\left( \mathbf{W}\right) }\left( s,t\right) ^{2}dsdt\right) %
\right] =\left( C_{\text{odd}}\left( 2u\right) \frac{4\mathcal{T}\left(
u\right) }{u}\right) ^{-\frac{1}{2}};$

\item (ii) $\mathbb{E}\left[ \exp \left( -\frac{u^{2}}{2}\int_{\left[ 0,1%
\right] ^{2}}\mathbf{B}_{0}^{\left( \mathbf{W}\right) }\left( s,t\right)
^{2}dsdt\right) \right] =\{S\left( u\right) \}^{-\frac{1}{2}};$

\item (iii) $\mathbb{E}\left[ \exp \left( -\frac{u^{2}}{2}\int_{\left[ 0,1%
\right] ^{2}}\mathbf{K}^{\left( \mathbf{W},1\right) }\left( s,t\right)
^{2}dsdt\right) \right] =\{S_{\text{odd}}\left( 2u\right) \}^{-\frac{1}{2}}.$
\end{enumerate}
\end{proposition}

\bigskip

Then, we have

\bigskip

\begin{theorem}
Under the above assumptions and notation, for every $u\in \Re $

\begin{enumerate}
\item
\begin{eqnarray}
&&\mathbb{E}\left[ \exp \left( iu\mathbf{I}\right) \right]  \label{Four1} \\
&=&\mathbb{E}\left[ \exp \left( -\frac{u^{2}}{2}\int_{\left[ 0,1\right] ^{2}}%
\left[ \mathbf{B}_{0}^{\left( \mathbf{W}\right) }\left( t_{1},t_{2}\right)
-\int_{\left[ 0,1\right] ^{2}}\mathbf{B}_{0}^{\left( \mathbf{W}\right)
}\left( u_{1},u_{2}\right) du_{1}du_{2}\right] ^{2}dt_{1}dt_{2}\right) %
\right]  \notag \\
&=&\left\{ C_{\text{odd}}\left( \frac{u}{2}\right) \frac{16\mathcal{T}\left(
u/4\right) }{u}\times S\left( \frac{u}{4}\right) \right\} ^{-\frac{1}{2}%
}\times S_{\text{odd}}\left( \frac{u}{2}\right) ,  \notag
\end{eqnarray}

\item
\begin{eqnarray}
&&\mathbb{E}\left[ \exp \left( iu\mathbf{J}\right) \right]  \label{Fouri2} \\
&=&\mathbb{E}\left[ \exp \left( -\frac{u^{2}}{2}\int_{\left[ 0,1\right] ^{2}}%
\left[ \mathbf{B}_{0}^{\left( \mathbf{W}\right) }\left( t_{1},t_{2}\right)
-\int_{\left[ 0,1\right] ^{2}}\mathbf{B}_{0}^{\left( \mathbf{W}\right)
}\left( t_{1},u_{2}\right) du_{2}\right] ^{2}dt_{1}dt_{2}\right) \right]
\notag \\
&=&\left\{ S\left( \frac{u}{2}\right) \right\} ^{-1}  \notag
\end{eqnarray}

\item
\begin{eqnarray}
&&\mathbb{E}\left[ \exp \left( iu\mathbf{Y}\right) \right]  \label{Fouri3} \\
&=&\mathbb{E}\left\{ \exp \left( \int_{\left[ 0,1\right] ^{2}}\left[ \mathbf{%
B}_{0}^{\left( \mathbf{W}_{1}\right) }\left( t_{1},t_{2}\right) -\int_{0}^{1}%
\mathbf{B}_{0}^{\left( \mathbf{W}_{1}\right) }\left( t_{1},u_{2}\right)
du_{2}\right. \right. \right.  \notag \\
&&\left. \left. \left. -\int_{0}^{1}\mathbf{B}_{0}^{\left( \mathbf{W}%
_{1}\right) }\left( u_{1},t_{2}\right) du_{1}+\int_{\left[ 0,1\right] ^{2}}%
\mathbf{B}_{0}^{\left( \mathbf{W}_{1}\right) }\left( u_{1},u_{2}\right)
du_{1}du_{2}\right] ^{2}dt_{1}dt_{2}\right) \right\}  \notag \\
&=&\left\{ S\left( \frac{u}{4}\right) \right\} ^{-2}  \notag
\end{eqnarray}
\end{enumerate}
\end{theorem}

\begin{proof}
The first equality in (\ref{Four1}) follows from conditioning and
independence. To obtain the second just recall that Theorem 3 implies that
\begin{eqnarray*}
&&\mathbb{E}\left[ \exp \left( -\frac{u^{2}}{2}\int_{\left[ 0,1\right] ^{2}}%
\left[ \mathbf{B}_{0}^{\left( \mathbf{W}\right) }\left( t_{1},t_{2}\right)
-\int_{\left[ 0,1\right] ^{2}}\mathbf{B}_{0}^{\left( \mathbf{W}\right)
}\left( u_{1},u_{2}\right) du_{1}du_{2}\right] ^{2}dt_{1}dt_{2}\right) %
\right] \\
&=&\mathbb{E}\left[ \exp \left( -\frac{\left( u/4\right) ^{2}}{2}\int_{\left[
0,1\right] ^{2}}\mathbf{B}_{0}^{\left( \mathbf{W}\right) }\left(
t_{1},t_{2}\right) ^{2}dt_{1}dt_{2}\right) \right] \times \mathbb{E}\left[
\exp \left( -\frac{\left( u/4\right) ^{2}}{2}\int_{\left[ 0,1\right] ^{2}}%
\mathbf{B}^{\left( \mathbf{W}\right) }\left( t_{1},t_{2}\right)
^{2}dt_{1}dt_{2}\right) \right] \\
&&\times \mathbb{E}\left[ \exp \left( -\frac{\left( u/4\right) ^{2}}{2}\int_{%
\left[ 0,1\right] ^{2}}\mathbf{K}^{\left( \mathbf{W},1\right) }\left(
t_{1},t_{2}\right) ^{2}dt_{1}dt_{2}\right) \right] ^{2},
\end{eqnarray*}%
and the conclusion follows from Proposition 5. Likewise,
\begin{eqnarray*}
&&\mathbb{E}\left[ \exp \left( -\frac{u^{2}}{2}\int_{\left[ 0,1\right] ^{2}}%
\left[ \mathbf{B}_{0}^{\left( \mathbf{W}\right) }\left( t_{1},t_{2}\right)
-\int_{\left[ 0,1\right] ^{2}}\mathbf{B}_{0}^{\left( \mathbf{W}\right)
}\left( t_{1},u_{2}\right) du_{2}\right] ^{2}dt_{1}dt_{2}\right) \right] \\
&=&\mathbb{E}\left[ \exp \left( -\frac{\left( u/2\right) ^{2}}{2}\int_{\left[
0,1\right] ^{2}}\mathbf{B}_{0}^{\left( \mathbf{W}\right) }\left(
t_{1},t_{2}\right) ^{2}dt_{1}dt_{2}\right) \right] ^{2},
\end{eqnarray*}%
so that the proof is achieved with another application of Proposition 5.
Formula (\ref{Fouri3}) is proved in exactly the same way.
\end{proof}

\bigskip

As pointed out in the Introduction, Theorem 6 extends part of the results
contained in \cite[Section 4]{DPY} and \cite{JuliaNualart}.

\bigskip

\textbf{Acknowledgements -- }The authors thank J.R.\ Pycke for showing them
the paper \cite{Pycke} prior to publication, and for several stimulating
discussions.

\end{document}